# Characterization of coincidence site lattices of oblique planar lattices

**Marco Antonio Rodríguez-Andrade**[(1)], **Gerardo Aragón-González**[(2*)], **José Luis Aragón Vera**[(3)]

[(1)]Departamento de Matemáticas, ESFM, IPN
[(2)]PDPA. UAM-Azcapotzalco.
[(3)]Centro de Física Aplicada y Tecnología Avanzada, Universidad Nacional Autónoma de México.
Apartado Postal 1-1010, 76000 Querétaro, México.
[*]E-mail: gag@correo.azc.uam.mx.

**Abstract**

Coincidence site lattices of oblique planar lattices are algebraically characterized using as basic tool the Cartan-Dieudonné theorem, that is, the decomposition of an orthogonal transformation as a product of reflections. The case of rectangular lattices is worked out in detail. We use the rectangular lattices for obtain the characterization of the corresponding obliques.

**1. Preliminaries**

A *lattice* $\Gamma$ in $\mathbf{R}^n$, consists of all integer linear combinations of a basis $\{\mathbf{a_1}, \mathbf{a_2}, \ldots, \mathbf{a_n}\}$ of $\mathbf{R}^n$, that is, $\Gamma = Z^n = Z\mathbf{a}_1 \oplus Z\mathbf{a}_2 \oplus \cdots \oplus Z\mathbf{a}_n$. In terms of group structure, leaving aside its geometric meaning, a lattice $\Gamma$ is a finitely-generated free Abelian group. In this sense, a subset $\Gamma' \subset \Gamma$ is called *sublattice* if it is a subgroup of finite index, where the index is the number of distinct cosets of $\Gamma'$ in $\Gamma$, denoted by $[\Gamma : \Gamma']$. The following two general definitions are central for the coincidence problem [1]:

**Definition 1.** *Two lattices $\Gamma_1$ and $\Gamma_2$ are called commensurate, denoted by $\Gamma_1 \sim \Gamma_2$, if and only if $\Gamma_1 \cap \Gamma_2$ is a sublattice of both $\Gamma_1$ and $\Gamma_2$.*

**Definition 2.** *Let $\Gamma$ be a lattice in $\mathbf{R}^n$. An orthogonal transformation $T \in O(n)$ is called a coincidence isometry of $\Gamma$ if and only if $T\Gamma \sim \Gamma$. The integer $\Sigma(T) := [\Gamma : \Gamma \cap T\Gamma]$ is called the coincidence index of $T$ with respect to $\Gamma$. If $T$ is not a coincidence isometry then $\Sigma(T) := \infty$. Two useful sets are also defined:*

$$OC(\Gamma) := \{T \in O(n) \mid \Sigma(T) < \infty\}, \quad SOC(\Gamma) := \{T \in OC(\Gamma) \mid \det(T) = 1\}.$$

Since we will be dealing with quotient groups, some words about the order of these groups should be said. The elements of the quotient group $T(\Gamma)/(\Gamma \cap T(\Gamma))$ will be denoted by $[a]$, where $a \in T(\Gamma)$, that is, $[a] = \{b \in T(\Gamma) \mid (a-b) \in \Gamma \cap T(\Gamma)\}$. In this case we have a finite Abelian group with identity $[0] = \Gamma \cap T(\Gamma)$. Since all the elements of this group have finite order, given $x \in \Gamma$, the order of $[T(x)]$, denoted as $o([T(x)])$, is finite and defined as

$$o([T(x)]) = \min\{k \in N \mid kT(x) \in \Gamma \cap T(\Gamma)\} \qquad (1)$$

In Ref. [2], the general problem of coincidence lattices in several planar lattices was worked out using Clifford algebras. In what follows we summarize the main results (without requiring Clifford algebras) that will be used in this work.

Suppose that $\mathbf{a} \in \mathbf{R}^n$ is a non-zero vector. By $\varphi_\mathbf{a}$ we denote a simple reflection (a reflection by a hyperplane orthogonal to $\mathbf{a}$). Thus, if $H_\mathbf{a}$ denotes the orthogonal complement to $\mathbf{a}$, then:

$$\varphi_\mathbf{a}(\mathbf{a}) = -\mathbf{a}; \quad \varphi_\mathbf{a}(\mathbf{w}) = \mathbf{w} \text{ if } \mathbf{w} \in H_\mathbf{a}$$

$$\varphi_\mathbf{a}(\mathbf{x}) = \mathbf{x} - 2\frac{\mathbf{x} \cdot \mathbf{a}}{\|\mathbf{a}\|^2}\mathbf{a} \tag{2}$$

By the Cartan Theorem[3], any given orthogonal transformation in $\mathbf{R}^n$ can be written as a product of simple reflections [4]. In our particular case, $\mathbf{a} \in \mathbf{R}^2$ and the hyperplane is a line perpendicular to $\mathbf{a}$. Also, we can always find $\mathbf{x}, \mathbf{y} \in \mathbf{R}^2$ such that $T_\theta = \varphi_\mathbf{x}\varphi_\mathbf{y}$.

In [2, Proposition 7] it was proved that for any lattice $\Gamma$ in $\mathbf{R}^2$ if $\varphi_\mathbf{c} \in OC(\Gamma)$ is a reflection then there exists a scalar $\lambda \neq 0$ such that $\lambda\mathbf{c} \in \Gamma$. Since $\varphi_\mathbf{c} = \varphi_{\lambda\mathbf{c}}$, this condition can be rephrased as follows:

**Lemma 1.** *A simple reflection $\varphi \in OC(\Gamma)$ if and only if there exists a vector $\mathbf{c} \in \Gamma$ such that $\varphi = \varphi_\mathbf{c}$.*

Now, for a rotation $T_\theta$ (henceforth, rotation with angle $\theta$) it was shown in the same reference (Proposition 11) that for simple lattices generated by the basis $\{\mathbf{a}_1, \mathbf{a}_2\}$, the coincidence site lattice (CSL) problem for this case is solved by the following Lemma:

**Lemma 2.** $T_\theta \in SOC(\Gamma)$ *if and only if there exists a vector $\mathbf{c} \in \Gamma$ such that $T_\theta = \varphi_\mathbf{c}\varphi_{\mathbf{b}_2}$, where $\mathbf{b}_2 = \mathbf{e}_1$ (square) $\mathbf{b}_2 = \mathbf{a}_2$ (rectangular) and $\mathbf{b}_2 = \mathbf{a}_1 - \mathbf{a}_2$ (rhombic [5]).*

This result arose from the fact that in the rhombic and hexagonal lattice the vectors $\mathbf{d}_1 = \mathbf{a}_1 + \mathbf{a}_2$ and $\mathbf{d}_2 = \mathbf{a}_1 - \mathbf{a}_2$ are orthogonal and, in fact, $\{\mathbf{d}_1, \mathbf{d}_2\}$ defines a rectangular sublattice of the rhombic lattice. Notice the resemblance with the approach followed by Ranganathan [6].

## 2. The CSL problem for the oblique lattice

Let $\Gamma_O = \mathbf{Z}\mathbf{e}_1 \oplus \mathbf{Z}(\sigma\mathbf{a})$ be an oblique lattice with $\sigma \in \mathbf{R}$ and $\mathbf{a} = (\cos\omega)\mathbf{e}_1 + (\sin\omega)\mathbf{e}_2$. Define $L_O = \begin{pmatrix} 1 & \sigma\cos\omega \\ 0 & \sigma\sin\omega \end{pmatrix}$ and let $P_O = \begin{pmatrix} 1 & \sigma\cos\omega \\ \sigma\cos\omega & \sigma^2 \end{pmatrix}$ be the structure matrix with respect to the ordered basis $\mathbf{B}_O = \{\mathbf{e}_1, \sigma\mathbf{a}\}$. Clearly, the square case corresponds to $\sigma = 1$ and $\omega = \frac{\pi}{2}$; the rectangular case to $\omega = \frac{\pi}{2}$; the hexagonal case to $\sigma = \frac{1}{2}$ and $\omega = \frac{\pi}{6}$; and the rhombic case to $\sigma = 1$.

The next Theorem characterizes the coincidence orthogonal transformations:

**Theorem 1.** $T \in OC(\Gamma)$, $(\Gamma = \mathbf{Z}\mathbf{a}_1 + \mathbf{Z}\mathbf{a}_2)$ *if and only if $[T]_\mathbf{B} \in M_{2\times 2}(\mathbf{Q})$ with respect to the basis $\mathbf{B} = \{\mathbf{a}_1, \mathbf{a}_2\}$.*

*Proof. The condition is necessary. We now use the fact that there is exists an isomorphism between $\Gamma/(T(\Gamma) \cap \Gamma)$ and $T(\Gamma)/(T(\Gamma) \cap \Gamma)$ (see [1]). Let $T \in OC(\Gamma)$, the quotient group $T(\Gamma)/(T(\Gamma) \cap \Gamma) = \{[T(\mathbf{a})] | \mathbf{a} \in \Gamma\}$ is finite, where the equivalence class is given by: $[T(\mathbf{a})] = \{\mathbf{b} \in \Gamma | T(\mathbf{b}) - T(\mathbf{a}) \in T(\Gamma) \cap \Gamma\}$; and $[0] = T(\Gamma) \cap \Gamma$. Thus, the elements of this group have*

*finite order, that is, for each* $[T(\mathbf{a})] \in \Gamma/(T(\Gamma) \cap \Gamma)$ *there exists* $k_\mathbf{a} \in \mathbf{Z}$ *such that:* $k_\mathbf{a}[T(\mathbf{a})] = [0]$ *or* $[k_\mathbf{a} T(\mathbf{a})] = [0]$. *Thus, for each* $\mathbf{a} \in \Gamma$ *there exists* $k_\mathbf{a} T(\mathbf{a}) \in T(\Gamma) \cap \Gamma$, *that is,* $k_\mathbf{a} T(\mathbf{a}) \in \Gamma$. *Assume that*

$$T(\mathbf{a}_1) = \gamma_{11} \mathbf{a}_1 + \gamma_{21} \mathbf{a}_2,$$
$$T(\mathbf{a}_2) = \gamma_{12} \mathbf{a}_1 + \gamma_{22} \mathbf{a}_2,$$

hence,

$$[T]_\mathbf{B} = \begin{pmatrix} \gamma_{11} & \gamma_{12} \\ \gamma_{21} & \gamma_{22} \end{pmatrix},$$

*Then, there exist* $k_i \in \mathbf{Z}$, $i = 1, 2$, *such that:*

$$k_1 T(\mathbf{a}_1) = k_1 \gamma_{11} \mathbf{a}_1 + k_1 \gamma_{21} \mathbf{a}_2 \in \Gamma,$$
$$k_2 T(\mathbf{a}_2) = k_2 \gamma_{12} \mathbf{a}_1 + k_2 \gamma_{22} \mathbf{a}_2 \in \Gamma,$$

*and since* $k_i \in \mathbf{Z}$ *(* $i = 1, 2$ *) then* $\gamma_{ij} \in \mathbf{Q}$ *(* $i = 1, 2;\ j = 1, 2$ *). Therefore:*

$$[T]_\mathbf{B} = \begin{pmatrix} \gamma_{11} & \gamma_{12} \\ \gamma_{21} & \gamma_{22} \end{pmatrix} \in M_{2 \times 2}(\mathbf{Q}).$$

*The condition is sufficient. Indeed, let* $[T]_\mathbf{B} = \begin{pmatrix} \gamma_{11} & \gamma_{12} \\ \gamma_{21} & \gamma_{22} \end{pmatrix} \in M_{2 \times 2}(\mathbf{Q})$ *and*

$k_i = \min\{m \in \mathbf{Z}^+ \mid mT(a_i) \in \Gamma\}$. *From this it follows that*

$$j[T(\mathbf{a}_i)] \neq l[T(\mathbf{a}_i)] \text{ for } 0 < j, l < k_i \ (i = 1, 2)$$

*Now consider* $\mathbf{a} = \alpha \mathbf{a}_1 + \beta \mathbf{a}_2 \in \Gamma$, *by applying the division algorithm, we have that there exist integers* $q_1, q_2, r_1, r_2$ *such that:*

$$\alpha = k_1 q_1 + r_1; \beta = k_2 q_2 + r_2,$$

*where* $0 \leq r_i < k_i$. *Thus,*

$$[T(\mathbf{a})] = r_1[T(\mathbf{a}_1)] + r_2[T(\mathbf{a}_2)],$$

*and we obtain that the quotient group* $T(\Gamma)/(T(\Gamma) \cap \Gamma)$ *has at most* $k_1 k_2$ *elements. Consequently, the group* $\Gamma/(T(\Gamma) \cap \Gamma)$ *is finite and* $T \in OC(\Gamma)$.

Now, for the general case of the oblique lattice $\Gamma_O$, we are interested in the necessary and sufficient for $\varphi_\mathbf{c} \in OC(\Gamma_O)$ $(\mathbf{c} = \alpha_1 \mathbf{e}_1 + \alpha_2 \sigma \mathbf{a}, \alpha_1, \alpha_2 \in \mathbf{Z}, \sigma \neq 0 \in \mathbf{R})$. For this goal, it is enough to calculate the matrix $[\varphi_\mathbf{c}]_{\mathbf{B}_O}$ from:

$$\varphi_\mathbf{c}(\mathbf{e}_1) = \mathbf{e}_1 - 2\left(\frac{\mathbf{e}_1 \cdot \mathbf{c}}{\|\mathbf{c}\|^2}\right)\mathbf{c} = \mathbf{e}_1 - 2\left(\frac{\alpha_1 + \alpha_2 \sigma \cos \omega}{\alpha_1^2 + \sigma^2 \alpha_2^2 + 2\alpha_1 \alpha_2 \sigma \cos \omega}\right)\mathbf{c} \tag{3}$$

$$\varphi_\mathbf{c}(\sigma \mathbf{a}) = \sigma \mathbf{a} - 2\left(\frac{\sigma \mathbf{a} \cdot \mathbf{c}}{\|\mathbf{c}\|^2}\right)\mathbf{c} = \sigma \mathbf{a} - 2\left(\frac{(\sigma^2 \alpha_2 + \alpha_1 \sigma \cos \omega)}{\alpha_1^2 + \sigma^2 \alpha_2^2 + 2\alpha_1 \alpha_2 \sigma \cos \omega}\right)\mathbf{c} \tag{4}$$

From Theorem 1, $\varphi_\mathbf{c} \in OC(\Gamma_O)$ if and only if

$$\frac{\alpha_1 + \alpha_2 \sigma \cos\omega}{\alpha_1^2 + \sigma^2 \alpha_2^2 + 2\alpha_1\alpha_2\sigma\cos\omega} \in \mathbf{Q},$$

$$\frac{(\sigma^2\alpha_2 + \alpha_1\sigma\cos\omega)}{\alpha_1^2 + \sigma^2\alpha_2^2 + 2\alpha_1\alpha_2\sigma\cos\omega} \in \mathbf{Q}, \tag{5}$$

The group $OC(\Gamma)$ contains at least two elements $T_\theta(\theta = 0, \pi)$, which will be called trivial. Now, we will found necessary and sufficient conditions for $OC(\Gamma_O)$ containing non-trivial elements. Indeed, suppose that the group $OC(\Gamma_O)$ has non-trivial elements and let $T_\theta \in SOC(\Gamma_O)$. Then,

$$[T_\theta]_{\mathbf{B}_O} = L_O^{-1}\begin{pmatrix} \cos\theta & -\sin\theta \\ \sin\theta & \cos\theta \end{pmatrix} L_O = \begin{pmatrix} -\frac{\sin(\theta-\omega)}{\sin\omega} & -\sigma\frac{\sin\theta}{\sin\omega} \\ \frac{1}{\sigma}\frac{\sin\theta}{\sin\omega} & \frac{\sin(\theta+\omega)}{\sin\omega} \end{pmatrix} \in M_{2\times 2}(\mathbf{Q}),$$

so:

$$\frac{\sigma\frac{\sin\theta}{\sin\omega}}{\frac{1}{\sigma}\frac{\sin\theta}{\sin\omega}} = \sigma^2 \in \mathbf{Q},$$

From (5), we have:

$$\sigma\cos\omega \in Q,$$

Therefore, we can establish the next Theorem:

**Theorem 2.** *The group $OC(\Gamma_O)$ has non-trivial elements if and only if $\sigma^2, \sigma\cos\omega \in \mathbf{Q}$.*

**Corollary 1.** *The group $OC(\Gamma)$ has non-trivial elements if and only if $\sigma^2 \in \mathbf{Q}$, for a rectangular lattice ($\Gamma = \Gamma_\square$), or $\cos\omega \in \mathbf{Q}$ for a rhombic lattice ($\Gamma = \Gamma_\lozenge$). Moreover, $\varphi_\mathbf{c} \in OC(\Gamma_\square)$ or $\varphi_\mathbf{c} \in OC(\Gamma_\lozenge)$. for any $\mathbf{c} \in \Gamma_\square$ or $\Gamma_\lozenge$.*

The following Lemma is valid for arbitrary planar lattices:

**Lemma 1.** *Let $\Gamma = \mathbf{Z}\mathbf{a}_1 \oplus \mathbf{Z}\mathbf{a}_2$ be a lattice, and let $\Gamma_D$ be the lattice generated by the basis $\mathbf{D} = \{\mathbf{d}_1, \mathbf{d}_2\}$, where $\mathbf{d}_1 = \mathbf{a}_1 + \mathbf{a}_2$ and $\mathbf{d}_2 = \mathbf{a}_1 - \mathbf{a}_2$ are the diagonals. Then $\Gamma_D$ is a sublattice of $\Gamma$ and $[\Gamma : \Gamma_D] = 2$. From this, it follows that any vector $\mathbf{x} \in \Gamma$, can be written as:*

$$\mathbf{x} = \mathbf{d} + k\mathbf{a}_1 \tag{6}$$

*where $\mathbf{d} \in \Gamma_D$ and $k = 0, 1$.*

*Proof: Since the diagonals of a rectangle are congruent and the diagonals of a rhombus are perpendiculars, then in all the cases we have:*

$$[d_1] = \begin{pmatrix} 1 \\ 1 \end{pmatrix}, \quad [d_2] = \begin{pmatrix} 1 \\ -1 \end{pmatrix}$$

*and a $s\det\begin{pmatrix} 1 & 1 \\ 1 & -1 \end{pmatrix} = 2$, then $[\Gamma : \Gamma_D] = 2$. For the second part, consider $\mathbf{x} = \beta_1\mathbf{a}_1 + \beta_2\mathbf{a}_2$, if $\beta_1, \beta_2 \in \mathbf{Z}$. then,*

$$\mathbf{x} = \frac{\beta_1 + \beta_2}{2}\mathbf{d}_1 + \frac{\beta_1 - \beta_2}{2}\mathbf{d}_2 \tag{7}$$

*If $\beta_1, \beta_2$ have the same parity, then $\beta_1 + \beta_2$ and $\beta_1 - \beta_2$ are even, and $\mathbf{x} \in \Gamma$ (so (6) is satisfied for $k = 0$). If $\beta_1, \beta_2$ have different parity, then $\beta_1 + \beta_2$ and $\beta_1 - \beta_2$ are odd, and we can find integers $m, n$ such that $\beta_1 + \beta_2 = 2m + 1$ and $\beta_1 - \beta_2 = 2n + 1$. From (7) follows that $\mathbf{x} = m\mathbf{d}_1 + n\mathbf{d}_2 + \mathbf{d}_1 + \mathbf{d}_2$, thus $\mathbf{x}$ can be written as $\mathbf{x} = \mathbf{d} + \mathbf{a}_1$.*

**Corollary 2.** $OC(\Gamma) = OC(\Gamma_D)$

From this Corollary it is inferred that for characterizing the group of orthogonal transformations in rectangular lattices it is enough to characterize the group of orthogonal transformation in rhombic lattices, and the reciprocal is also true.

**Remark 1.** *Considering the rectangular lattice* $\Gamma_\Box$ *and let $\gamma$ be the angle between* $\mathbf{d}_1$ *and* $\sigma\mathbf{e}_1$. *From Corollary 2 it follows that* $\Gamma_{\Box D}$ *is a rhombic lattice and the angle between* $\mathbf{d}_1$ *and* $\mathbf{d}_2$ *is* $2\gamma$. *Even more,*

$$\cos\gamma = \frac{\sigma}{\sqrt{1+\sigma^2}}; \quad \tan 2\gamma = \frac{\sigma^2 - 1}{\sigma^2 + 1}.$$

*Now, if we consider to* $\Gamma_\Diamond$ *and to its respective* $\Gamma_{\Diamond D}$ *which is a rectangular lattice, we have:* $\dfrac{\|\mathbf{d}_1\|}{\|\mathbf{d}_2\|} = \sqrt{\dfrac{1+\cos\omega}{1-\cos\omega}}$.

In conclusion to characterize $OC(\Gamma_O)$ it is enough to study $OC(\Gamma_\Box)$. This is done in the following section.

## 3. Characterization of coincidence site lattices of rectangular planar lattices.

Let $\Gamma_\Box = \mathbf{Z}\mathbf{e}_1 + \mathbf{Z}(\sigma\mathbf{e}_2)$ be a rectangular lattice and $T_\theta$ be a rotation. By applying the algorithm proposed in [4] $T_\theta$ can be decomposed as a product of two simple reflections:

$$T_\theta = \varphi_\mathbf{c}\varphi_{\sigma\mathbf{e}_2} \tag{8}$$

where $\mathbf{c} = \alpha_1\mathbf{e}_1 + \alpha_2(\sigma\mathbf{e}_2) = \lambda(T(\mathbf{e}_1) - \mathbf{e}_1), \alpha_1, \alpha_2 \in \mathbf{Z}$; and $\lambda \neq 0 \in \mathbf{R}$. Since $\mathbf{B}_\Box = \{\mathbf{e}_1, \sigma\mathbf{e}_2\}$ is the basis of the lattice, from (3) and (4) we have:

$$[\varphi_\mathbf{c}]_{\mathbf{B}_\Box} = \begin{pmatrix} \frac{(-\alpha_1+\sigma\alpha_2)(\alpha_1+\sigma\alpha_2)}{\alpha_1^2+\sigma^2\alpha_2^2} & -2\sigma^2\frac{\alpha_1\alpha_2}{\alpha_1^2+\sigma^2\alpha_2^2} \\ -2\frac{\alpha_1\alpha_2}{\alpha_1^2+\sigma^2\alpha_2^2} & \frac{(\alpha_1-\sigma\alpha_2)(\alpha_1+\sigma\alpha_2)}{\alpha_1^2+\sigma^2\alpha_2^2} \end{pmatrix}, \quad [\varphi_{\sigma\mathbf{e}_2}]_{\mathbf{B}_\Box} = \begin{pmatrix} 1 & 0 \\ 0 & -1 \end{pmatrix}.$$

From (8) it follows that:

$$[T_\theta]_{\mathbf{B}_\Box} = \begin{pmatrix} \frac{(\alpha_1+\sigma\alpha_2)(-\alpha_1+\sigma\alpha_2)}{\alpha_1^2+\sigma^2\alpha_2^2} & 2\sigma^2\frac{\alpha_1\alpha_2}{\alpha_1^2+\sigma^2\alpha_2^2} \\ -2\frac{\alpha_1\alpha_2}{\alpha_1^2+\sigma^2\alpha_2^2} & \frac{(\alpha_1+\sigma\alpha_2)(-\alpha_1+\sigma\alpha_2)}{\alpha_1^2+\sigma^2\alpha_2^2} \end{pmatrix},$$

where the angle between the vectors $\mathbf{c}$ and $\sigma\mathbf{e}_2$ must be $\frac{\theta}{2}$. This implies that the tranformation $\varphi_\mathbf{c}\varphi_{\sigma\mathbf{e}_2}$ is $T_\theta$. As conclusion of all the above, we have:

**Theorem 3.** *Let* $\Gamma_\Box = \mathbf{Z}\mathbf{e}_1 + \mathbf{Z}(\sigma\mathbf{e}_2)$ *be a rectangular lattice with* $\sigma^2 \in \mathbf{Q}$. *Then* $T_\theta \in SOC(\Gamma_\Box)$ *if and onyl if there exist* $p, q \in \mathbf{Z}$ $(\gcd(p,q) = 1)$ *such that* $\dfrac{\theta}{2} = \arctan\left(\dfrac{p}{q}\dfrac{1}{\sigma}\right)$. *Furthermore,* $T_\theta = \varphi_\mathbf{c}\varphi_{\sigma\mathbf{e}_2}$ $(\mathbf{c} = p\mathbf{e}_1 + q\sigma\mathbf{e}_2)$ *and its matrix with respect to the basis* $\mathbf{B}_\Box$ *is given by,*

$$[T_\theta]_{\mathbf{B}_\Box} = \begin{pmatrix} \frac{(p+\sigma q)(-p+\sigma q)}{p^2+\sigma^2 q^2} & 2\sigma^2\frac{pq}{p^2+\sigma^2 q^2} \\ -2\frac{pq}{p^2+\sigma^2 q^2} & \frac{(p+\sigma q)(-p+\sigma q)}{p^2+\sigma^2 q^2} \end{pmatrix}$$

Now consider the oblique lattice $\Gamma_O = \mathbf{Z}\mathbf{e}_1 \oplus \mathbf{Z}(\sigma\mathbf{a}); \sigma \in \mathbf{R}$ and $\mathbf{a} = (\cos\omega)\mathbf{e}_1 + (\sin\omega)\mathbf{e}_2,$. In this case we have $\left[\varphi_{\mathbf{e}_2}\right]_{\mathbf{B}_O} = \begin{pmatrix} 1 & 2\sigma\cos\omega \\ 0 & -1 \end{pmatrix}$. From Theorem 2, $OC(\Gamma_O)$ has non-trivial elements if $\sigma^2, 2\sigma\cos\omega \in \mathbf{Q}$, thus $\varphi_{\mathbf{e}_2} \in OC(\Gamma_O)$, which implies that there exists $\lambda \in \mathbf{Q}$ such that $\lambda\mathbf{e}_2 \in \Gamma_O$, i.e. $\Gamma_O$ has a rectangular sublattice of finite index. Applying the algorithm given in [4], the rotation $T_\theta$ can be written as $T_\theta = \varphi_{\mathbf{c}} \varphi_{\lambda\mathbf{e}_2}$ and by the carried out above for rectangular case (proof of the Theorem 3), we can establish the following Lemma:

**Lemma 2.** $T_\theta \in SOC(\Gamma_O)$ if and only if there exists a vector $\mathbf{c} \in \Gamma_O$, such that $T_\theta = \varphi_{\mathbf{c}} \varphi_{\lambda\mathbf{e}_2}$. Even more: $[T_\theta]_{\mathbf{B}_0} = \begin{pmatrix} \frac{(\alpha_1+\sigma\alpha_2)(-\alpha_1+\sigma\alpha_2)}{\alpha_1^2+\sigma^2\alpha_2^2+2\sigma\alpha_1\alpha_2\cos\omega} & 2\sigma^2\alpha_2 \frac{\alpha_1+\sigma\alpha_2\cos\omega}{\alpha_1^2+\sigma^2\alpha_2^2+2\sigma\alpha_1\alpha_2\cos\omega} \\ -2\alpha_2 \frac{\alpha_1+\sigma\alpha_2\cos\omega}{\alpha_1^2+\sigma^2\alpha_2^2+2\sigma\alpha_1\alpha_2\cos\omega} & -\frac{(\alpha_1-\sigma\alpha_2+2\sigma\alpha_2\cos\omega)(\alpha_1+\sigma\alpha_2+2\sigma\alpha_2\cos\omega)}{\alpha_1^2+\sigma^2\alpha_2^2+2\sigma\alpha_1\alpha_2\cos\omega} \end{pmatrix}$

As in the rectangular case, from the Lemma 2 we obtain:

**Theorem 4.** $T_\theta \in SOC(\Gamma_O)$ if and only there exist $p,q \in \mathbf{Z}$ $(\gcd(p,q)=1)$ such that:

$$\left(\frac{\theta}{2}\right) = \arctan\left(\frac{p}{q}\tan\omega\right)$$

*Proof;* Let $\gamma$ the angle between $\mathbf{c} = \alpha_1\mathbf{e}_1 + \alpha_2(\sigma\mathbf{a})$ and $\mathbf{e}_2$. Since $\tan(\gamma) = \frac{\alpha_2\sigma\sin\omega}{\alpha_1+\alpha_2\sigma\cos\omega}$, the product $\varphi_{\mathbf{c}}\varphi_{\mathbf{e}_2}$ corresponds to a rotation by the angle $2\gamma$. A necessary and suffienct condition for $T_\theta \in SOC(\Gamma_O)$ is is that there exist $\alpha_1, \alpha_2 \in \mathbf{Z}$, such that:

$$\left(\frac{\theta}{2}\right) = \arctan\left(\frac{\alpha_2\sigma\cos\omega}{\alpha_1+\alpha_2\sigma\cos\omega}\tan\omega\right)$$

**4. Conclusions**

The problem of the coincidence site lattices of oblique planar lattices is addressed in this work. Our main result is that the necessary and sufficient conditions for an orthogonal transformation to be a a coincidence transformation can be found by studying the coincidence site lattices of rectangular planar lattices. We should finally mention that Theorem 5 is not applicable to the case of rhombic lattices since in that case $\omega = \frac{\pi}{2}$. It can however be worked out in a similar way as the rectangular lattice. Further work is underway.